\documentclass{amsart}
\usepackage{amsmath,amsthm,amssymb,IMjournal}

\usepackage{constants, color, mathtools, tikz, xcolor}
\usepackage{graphicx,hyperref}
\usepackage{inputenc} 
\usepackage{tikz}
\usepackage{rotating}
\usepackage{multirow} 
\usepackage{colortbl} 
\usetikzlibrary{shapes,snakes}
\usepackage{hyperref}
\usepackage[linesnumbered,ruled,vlined]{algorithm2e}
\SetKwInput{KwInput}{Input}                
\SetKwInput{KwOutput}{Output}              

\begin{document}
\newtheorem{thm}{Theorem}[section]
\newtheorem{rem}{Remark}[section]
\newtheorem{prop}{Proposition}[section]
\newtheorem{conj}{Conjecture}[section]

\numberwithin{equation}{section}

\title{Sparse Sets in Triangle-free Graphs}

\author{\|T\i naz |Ekim|, Istanbul,
        \|Burak N. |Erdem|, Istanbul,
        \|John |Gimbel|, Fairbanks}

\vspace*{-2.5em}
\begin{center}
  \small June, 2025
\end{center}
\vspace{1em}


\abstract 
A set of vertices is $k$-sparse if it induces a graph with a maximum degree of at most $k$.  In this missive, we consider the order of the largest $k$-sparse set in a triangle-free graph of fixed order.  We show, for example, that every triangle-free graph of order 11 contains a 1-sparse 5-set; every triangle-free graph of order 13 contains a 2-sparse 7-set; and every triangle-free graph of order 8 contains a 3-sparse 6-set.  Further, these are all best possible.  

For fixed $k$, we consider the growth rate of the largest $k$-sparse set of a triangle-free graph of order $n$.  Also, we consider Ramsey numbers of the following type.  Given $i$, what is the smallest $n$ having the property that all triangle-free graphs of order $n$ contain a 4-cycle or a $k$-sparse set of order $i$. We use both direct proof techniques and an efficient graph enumeration algorithm to obtain several values for defective Ramsey numbers and a parameter related to largest sparse sets in triangle-free graphs, along with their extremal graphs.
\endabstract



\thanks
   First and second authors have been supported by the Scientific and Technological Research Council of Turkey TUBITAK Grant No 122M452. This joint work was initiated when John Gimbel visited the first and second authors with the support of the Turkish Mathematical Society (TMD) and the Istanbul Center for Mathematical Sciences (IMBM). The support of both institutions is greatly acknowledged. John Gimbel also received funding from the European Research Council (ERC) under the European Union's Horizon 2020 research and innovation programme (grant agreement No 810115).
\endthanks

\section{Introduction}
Given positive integers $i$ and $j$, the Ramsey number $R(i, j)$ is the smallest
natural number such that every graph of order at least $R(i, j)$ has a clique of order $i$ or an independent set of order $j$. These so called classical Ramsey numbers along with a number of variations are extensively
studied in the literature. Among various generalizations 
we find so called defective Ramsey numbers that have been the focus of several research papers \cite{AkdemirEkim, gimbelson, 1-dependent, ExactDefRamClasses, SubExt, ekim2013some, SmallDefRamPerfect}. This variation relaxes the notions of cliques and independent sets in the following way. A \textit{$k$-sparse set} is a set of vertices that induces a graph with maximum degree $k$ or less.  A \textit{$k$-dense set} is the complement of a $k$-sparse set.  In other words, each vertex in a $k$-dense set “misses” at most $k$ other vertices in its neighborhood.  A $k$-sparse ($k$-dense) \textit{$j$-set} is a $k$-sparse ($k$-dense) set of order $j$. A set is \textit{$k$-defective} (or \textit{$k$-uniform}) if it is a $k$-sparse or $k$-dense set. The \textit{defective Ramsey number} $R^\mathcal{G}_k (i, j)$ is the least $n$ such that all graphs of order $n$ in the class $\mathcal{G}$
have either a $k$-dense $i$-set or a $k$-sparse $j$-set.

In \cite{1-dependent}, some specific 1-defective Ramsey numbers are derived under a somewhat different terminology. In \cite{gimbelson, ekim2013some}, additional 1-defective and 2-defective Ramsey numbers are found
using direct proofs. Further, several bounds are displayed  for defective Ramsey numbers. It seems direct proofs have reached their limits in finding
new values of defective Ramsey numbers. Indeed, this is rather not surprising given the
great difficulty in computing specific defective Ramsey numbers. Having observed this
fact, some computer based generation methods are used in \cite{AkdemirEkim, gimbelson, SubExt} to improve 
known bounds on defective Ramsey numbers (and certain other defective parameters).

Noting that computing exact Ramsey numbers is extremely unlikely, various approaches
are adopted in mathematical literature to partially deal with this problem. One tactic is to consider restricted graph families. In \cite{RamseyPlanar} for example, all classical Ramsey numbers in planar graphs
are found. In \cite{RamseyDegree1, RamseyDegree2, RamseyDegree3}, the authors compute several
Ramsey numbers for graphs with bounded degree. Ramsey numbers for claw-free graphs are discussed in \cite{RamseyClawfree}. After all these studies dating back to 1980's and 1990's, this approach
seems to have become popular again. In \cite{ramseygraphclasses}, we find a systematic study of Ramsey numbers in various graph classes. It seems that computation for claw-free graphs is as difficult as it is for arbitrary graphs. Further, \cite{ramseygraphclasses} exhibits all classical
Ramsey numbers for perfect graphs and some well-known subclasses of claw-free graphs.
We note another work \cite{defcolperfect} focuses on the complexity of the coloring problem where
each color class is a $k$-sparse set (called the $k$-defective coloring problem) when restricted to subclasses of perfect graphs.

Recently, the approach of considering Ramsey numbers in various graph classes has been applied to defective Ramsey numbers. In \cite{SmallDefRamPerfect}, Ekim et al. present small 1-defective Ramsey numbers for perfect graphs. In \cite{ExactDefRamClasses}, Demirci et al. study $k$-defective Ramsey numbers (for any $k$) and provide exact formulas for forests, cacti, bipartite graphs, split graphs and cographs. They provide conjectures for the few exceptions left as open questions. In both of these studies \cite{SmallDefRamPerfect} and \cite{ExactDefRamClasses}, the authors observe that the limits of direct proof techniques seem to be reached. As such, Demirci et al. focus more recently in \cite{SubExt} on the computation of defective Ramsey numbers by combining their efficient graph generation algorithm, called Sub-extremal, with classical proof techniques. They provide new defective Ramsey numbers in perfect graphs, bipartite graphs, and chordal graphs.

In this paper, we investigate defective Ramsey numbers in triangle-free graphs. Our contributions are two-fold: we provide both direct proofs and computer assisted results using an efficient implementation of the Algorithm Sub-extremal from \cite{SubExt} adapted for triangle-free graphs. Some simple observations show that dense sets in triangle-free graphs are very restricted. This implies that we rather focus on the existence of sparse sets of given size. From \cite{ajtai} we know that every triangle-free graph of order $n$ contains an independent set of order at least $\sqrt{n \log n}$ and by \cite{kim} this is asymptotically best possible.  In this work, we extend this notion and consider large sparse sets in triangle-free graphs.  We work out some specific values and produce computer assisted proofs of others. All defective Ramsey numbers in triangle-free graphs obtained with direct proof techniques in Sections \ref{sec:T} and \ref{sec:R}, as well as those obtained using an efficient graph generation approach in Section \ref{sec:computer} are summarized in Tables \ref{table:T_k(j)}, \ref{table:R_{1}(i,j)}, \ref{table:R_{2}(i,j)}, \ref{table:R_{3}(i,j)}, \ref{table:R_{4}(i,j)}. We postpone a more detailed description of our results until the end of Section \ref{sec:prem}, after all formal definitions, notations and preliminary remarks are introduced.

\section{Definitions and preliminary remarks}\label{sec:prem}

Let $G=(V,E)$ be a graph. We will denote the order $|V|$ of a graph by $n$. A subgraph $H\subseteq G$ is a graph on $V' \subseteq V$ and $E' \subseteq E$ with both end-vertices of each edge of $E'$ in $V'$. If all edges with both end-vertices in $V'$ are in $E'$, then $H$ is said to be an \textit{induced} subgraph of $G$. In our context, whenever we say that a graph contains a subgraph, we  always mean as a partial subgraph, unless stated otherwise. For a vertex $x\in V$, we denote by $N(x)$ the set of \textit{neighbors} of $x$, that is, vertices adjacent to $x$. The \textit{degree} of a vertex $x$ is $d(x)=|N(x)|$. The maximum degree of a graph, denoted by $\Delta(G)$, is the largest vertex degree in $G$. Likewise, the minimum degree of a graph, denoted by $\delta(G)$, is the smallest vertex degree in $G$. We also have $N[x]=N(x)\cup \{x\}$. For a vertex $x\in V$ and a subgraph $H\subseteq G$, we denote by $N_H(x)$ the set of neighbors of $x$ in $H$, that is $N(x)\cap V(H)$. Similarly, the degree of $x$ in $H$ is $d_H(x)=|N_H(x)|$. For a subset of vertices $X\subset V$, the \textit{neighborhood} of $X$, denoted by $N(X)$, is defined as $N(X)=(\cup_{x\in X}N(x))\setminus X$. 

For a graph $G$ and a subgraph $H$, we use the notation $G\setminus H$ to mean the subgraph of $G$ induced by all vertices in $V(G)\setminus V(H)$. We also use the same notation when we remove a set of vertices from a graph. For  graphs $H$ and $G$, we say that $G$ is \textit{$H$-free} if it does not contain $H$ as an induced subgraph. A \textit{path} on $n$ vertices is denoted by $P_n$, and a \textit{cycle} on $n$ vertices, also called an \textit{$n$-cycle}, is denoted by $C_n$. A complete bipartite graph on $p$ and $\ell$ vertices in each part is denoted by $K_{p,\ell}$. The \textit{distance} between two vertices is the length of a shortest path between them. The \textit{girth} of a graph $G$, denoted by $g(G)$, is the length of a shortest induced cycle in it. A set of vertices is called \textit{independent} if all vertices in it are pairwise non-adjacent. Generalizing the notation for the size of a largest independent set $\alpha(G)$, we adopt the notation $\alpha_k(G)$ to denote the size of a largest $k$-sparse set of a graph $G$. An \textit{extremal} graph for $R^{\mathcal{G}}_k(i,j)$ is a graph in the class $\mathcal{G}$ on $R^{\mathcal{G}}_k(i,j)-1$ vertices containing neither a $k$-dense $i$-set nor a $k$-sparse $j$-set. We use the notation $R^{\Delta}_k(i,j)$ to denote the $k$-defective Ramsey numbers in triangle-free graphs.

In this work, sparse sets in triangle-free graphs will be our main focus. This is justified by the following remarks. First, let us note that in general, finding such sets is difficult; we know from \cite{ekim2013some} that finding a largest $k$-sparse set for any fixed $k$ is $NP$-complete even in restricted cases:

\begin{thm}\label{thm:NPc}
\cite{ekim2013some}
For fixed $k\geq 2$, given a graph $G$ and an integer $t$, the problem of deciding if $\alpha_k(G)\geq t$ is NP-complete. The result holds when restricted to planar graphs with maximum degree $k+1$ and girth $g$, where $g$ is arbitrarily large.
\end{thm}

Moreover, Theorem \ref{thm:NPc} holds for $k=1$ when the maximum degree is three \cite{ekim2013some}. For $k=0$, deciding if there is an independent set of size at least $t$ is NP-complete in several restricted cases including triangle-free graphs \cite{poljak}.

We proceed with some observations on the absence of large dense sets in triangle-free graphs. This motivates the study of sparse sets in triangle-free graphs in further sections.

\begin{rem}
In a triangle-free graph $G$, a 1-dense 4-set can only be a $C_4$. Morover, $G$ does not admit 1-dense $i$-set for any $i\geq 5$. 
\end{rem}

In a similar way, we can show that the only 2-dense 5-sets are $C_5$ and $K_{2,3}$, and the only 2-dense 6-set is $K_{3,3}$. These observations can be generalized for $k\geq 2$ as follows.

\begin{prop}\label{prop:Tnotation}
In a triangle-free graph, there is no $k$-dense $i$-set for $i\geq 2k+3$. Moreover, this bound is best possible and the unique $k$-dense triangle-free graph on $2k+2$ vertices is $K_{k+1,k+1}$.
\end{prop}
\begin{proof}
Let $G$ be a triangle-free graph and assume to the contrary it contains a $k$-dense set, say $A$, having at least $2k+3$ vertices. Let $x$ be a vertex of $A$. Note $x$ can miss at most $k$ other vertices of $A$. Hence $x$ is adjacent to at least $k+2$ other vertices of $A$. Since $G$ is triangle-free, $N(x)\cap A$ is an independent set. But then, a vertex in $N(x)\cap A$ misses the other $k+1$ vertices of $N(x)\cap A$, contradicting the fact that $A$ is $k$-dense. 

Note, a $k$-dense graph of order $2k+2$ contains at least $(k+1)^2$ edges.  By Turan's Theorem \cite{turan}, there is only one triangle-free graph of order $2k+2$ on $(k+1)^2$ edges, namely $K_{k+1,k+1}$.  Further, if a graph of order $2k+2$ contains more than $(k+1)^2$ edges, it must contain a triangle.  Hence, our result is best possible.
\end{proof}

Since there is no $k$-dense $i$-set for $i\geq 2k+3$ in a triangle-free graph, Proposition \ref{prop:Tnotation} implies that for each $k$, we have $R^{\Delta}_k(i,j)=R^{\Delta}_k(i',j)$ for all $i,i' \geq 2k+3$. Without $k$-dense sets, it makes sense to focus on $k$-sparse sets. This suggests the following notation. Let $T_k(j)$ be the minimum order $n$ such that every triangle-free graph of order $n$ has a $k$-sparse set of size $j$. We would say an \textit{extremal} graph for $T_k(j)$ is a triangle-free graph with $T_k(j)-1$ vertices having no $k$-sparse set of order $j$. With this notation, we have $R^{\Delta}_k(i,j)=T_k(j)$ for all $i\geq 2k+3$. Motivated by this, we proceed by proving some exact values for $T_k(j)$ in Section \ref{sec:T}. Then, in Section \ref{sec:R}, we focus on $R_1^{\Delta}(3,j)$ for $j\geq 3$, and $R_1^{\Delta}(4,j)$ for $j\geq 4$, the only 1-defective Ramsey numbers of interest for triangle-free graphs. Both Sections \ref{sec:T} and \ref{sec:R} contain results shown by classical proof techniques. In Section \ref{sec:computer}, we compute several new values by efficient computer enumeration techniques. Based on these results, we conjecture that $T_k(k+i)=k+2i-1$ for all $i$ and $k$ such that $2\leq i \leq k$. All of our codes and the extremal graphs we obtain are available online at \cite{github}.

\section{Sparse sets in triangle-free graphs}\label{sec:T}

The following lower bound on $\alpha_k(G)$ allows us to derive some values of $T_k(n)$. Note that the following lower bound is for general graphs, not restricted to triangle-free graphs.

\begin{prop}\label{prop:alphaFormula} 
For a graph $G$ and fixed $k$, we have $\alpha_k(G)\geq \frac{n}{\bigl\lceil \frac{\Delta(G)+1}{k+1}\bigr\rceil}$.
\end{prop}
\begin{proof}
We rely on a proof technique found in \cite{lovasz}. Set $j=\bigl\lceil \frac{\Delta(G)+1}{k+1}\bigr\rceil$. Color the vertices of $G$ with $j$ colors so that the number of monochromatic edges (those edges having the same color on both end-vertices) is minimized. We claim that every color class is a $k$-sparse set. Assume this does not hold, that is, there is a vertex $x$ with at least $k+1$ neighbors of the same color as $x$. Then one of the remaining $j-1$ colors, say $c$, occurs at most $k$ times in the neighborhood of $x$, since otherwise $d(x) \geq \Delta(G)+1$, a contradiction. By recoloring $x$ with $c$ we obtain a coloring of $G$ with fewer monochromatic edges, a contradiction. Now, by the Pigeonhole Principle, one of the color classes has at least $\frac{n}{\bigl\lceil \frac{\Delta(G)+1}{k+1}\bigr\rceil}$ vertices, and the proof is complete.
\end{proof}

Before investigating specific values of $T_k(n)$, it is worth noting the case where $k = 0$, in which, the sparse set in question is an independent set. Consequently, the parameter $T_0(j)$ is equivalent to the classical Ramsey number $R(3,j)$. Moving onto the sparse sets, the first non-trivial value is $T_1(3)=5$. Note $C_4$ is an extremal graph on 4 vertices. Suppose $G$ is a triangle-free graph of order 5. If $G$ is bipartite, it contains an independent set on 3 vertices. So, suppose $G$ is not bipartite. Note $G$ contains an odd cycle which is not a triangle. Hence $G$ contains a 5-cycle without chord. This graph contains a 1-sparse 3-set. 

In the following, we will repeatedly use (without explicitly mentioning it) the observation that any open neighborhood in a triangle-free graph is independent. 

\begin{thm}\label{thm:T1_4_7} 
With the preceeding notation, $T_1(4)=7$.
\end{thm}
\begin{proof}
$K_{3,3}$ is a triangle-free graph of order 6 which does not contain a 1-sparse 4-set. Thus, $T_1(4)\geq 7$.
Let $G$ be a triangle-free graph of order 7. If ${\Delta(G)}\geq 4 $, then $N(x)$ contains a 1-sparse 4-set. If ${\Delta(G)}\leq 3 $, by Proposition \ref{prop:alphaFormula}, the cardinality of a 1-sparse set is at least $4$. Consequently, every triangle-free graph with $7$ vertices includes a 1-sparse 4-set. 
Hence, the desired result.
\end{proof}

\begin{thm}
With the preceeding notation, $T_1(5)=11$. 
\end{thm}
\begin{proof}
The blow-up of a $C_5$ where every vertex is replaced with two independent vertices is a graph or order 10 which contains no 1-sparse 5-set. Thus, $T_1(5) \geq 11$.

Let $G$ be a triangle-free graph of order 11. If $G$ contains a vertex $x$ of degree at least 5, then $N(x)$ contains a 1-sparse 5-set. So suppose $\Delta(G)\leq 4$. If $G$ has a vertex $x$ of degree 3, then $V\setminus N[x]$ has 7 vertices and contains a 1-sparse 4-set, say $A$, by Theorem \ref{thm:T1_4_7}. Now, $A\cup \{x\}$ is a 1-sparse 5-set. So, suppose $G$ has no vertex of degree three and similarly no vertex of degree less than three.
 
So, assume $G$ is a 4-regular triangle-free graph of order 11. For some vertex $x$, let $A=N(x)$ and $B=V\setminus N[x]$. Note $A$ is independent. We have $|A|=4$, $|B|=6$, where each vertex of $A$ is adjacent with exactly three vertices in $B$. Thus, there are exactly 12 edges with one end-vertex in $A$ and the other in $B$. If there is a vertex $b\in B$ having at most one neighbor in $A$, then $A\cup \{b\}$ is a 1-sparse 5-set. Otherwise, every vertex of $B$ has exactly 2 neighbors in $A$. This implies that the graph induced by $B$ is 2-regular; thus a 6-cycle (since triangles are forbidden). Taking a 1-sparse 4-set in this 6-cycle together with $x$ yields a 1-sparse 5-set. Hence, the desired result.
\end{proof}

In the sequel, we study 2-sparse sets in triangle-free graphs. We start with the fırst non-trivial value of the $T_k(j)$ for all $k\geq 2$.

\begin{thm}\label{thm:T_k(k+2)}
With the preceeding notation, $T_k(k+2)=k+3$, for $k \geq 2$ with $K_{1,k+1}$ as the unique extremal graph. 
\end{thm}
\begin{proof}
$K_{1,k+1}$ is a graph which do not include any triangles nor any $k$-sparse $(k+2)$-sets. Therefore, $T_k(k+2) \geq k+3$.

Let $G$ be a triangle-free graph of order $k+3$ and let $x$ be a vertex of maximum degree. If $\Delta(G) \leq k$ then $G$ is $k$-sparse. If $\Delta(G) \geq k+2$, then $N(x)$ is independent and hence contains a $k$-sparse $(k+2)$-set. Suppose $\Delta(G) = k+1$, and let $y$ be the vertex in $V \setminus N[x]$. Since $k\geq 2$, the set $\{x,y\}$ together with any $k$ vertices from $N(x)$ is a $k$-sparse $(k+2)$-set. Thus, $T_k(k+2)=k+3$.

Let us now show that $K_{1,k+1}$ is the unique extremal graph. Indeed, by the previous observation, an extremal graph with $k+2$ vertices has maximum degree $k+1$ or else it has a $k$-sparse $(k+2)$-set. Since it is a triangle-free graph, it can only be a $K_{1,k+1}$.
\end{proof}

Theorem \ref{thm:T_k(k+2)} implies in particular that $T_2(4)=5$. We proceed with the next values.
\begin{thm}
With the preceeding notation, $T_2(5)=9$. 
\end{thm}
\begin{proof}
$K_{4,4}$ contains no 2-sparse 5-set. Thus, $T_2(5)\geq 9$.
Let $G$ be a triangle-free graph of order 9. If ${\Delta(G)} \geq 5$, the neighborhood of a maximum degree vertex contains a 2-sparse 5-set. If ${\Delta(G)} \leq 4$, there exists a 2-sparse set of size at least 5, by Proposition \ref{prop:alphaFormula}. Hence, the desired result.
\end{proof}

\begin{thm}
With the preceeding notation, $T_2(6)=11$.
\end{thm}
\begin{proof}
$K_{5,5}$ contains no 2-sparse 6-set. Thus, $T_2(6)\geq 11$.
Let $G$ be a triangle-free graph of order 11. If ${\Delta(G)} \geq 6$, the neighborhood of a maximum degree vertex contains a 2-sparse 6-set. If ${\Delta(G)} \leq 5$, there exists a 2-sparse set of size at least 6, by Proposition \ref{prop:alphaFormula}. Hence, the desired result.
\end{proof}

\begin{thm}\label{thm:T_2(7)=13}
With the preceeding notation, $T_2(7)=13$.
\end{thm}
\begin{proof}
$K_{6,6}$ is a triangle-free graph of order 12 which has no 2-sparse 7-set. Thus ${T_2(7)\geq 13}$. 
Let $G$ be a triangle-free graph of order 13. If ${\Delta(G)} \geq 7$, the neighborhood of a maximum degree vertex contains a 2-sparse 7-set. If ${\Delta(G)} \leq 5$, there exists a 2-sparse set of size at least 7, by Proposition \ref{prop:alphaFormula}. So assume ${\Delta(G)} = 6$ and let $x$ be a vertex of degree six. Let $N(x)=\{v_1, v_2, v_3, v_4, v_5, v_6\}$ and $U=V\setminus N[x]=\{u_1, u_2, u_3, u_4, u_5, u_6\}$ as shown in Figure \ref{fig:proof_t_2_7}. If $U$ is a 2-sparse 6-set, then $\{x\} \cup U$ is a 2-sparse 7-set. Otherwise, there exists a vertex in $U$, say $u_1$, which is adjacent to at least three other vertices in $U$. Let $\{u_2, u_3, u_4\} \subseteq N(u_1)$. 
If $u_1$ has at most two neighbors in $N(x)$, then $ \{u_1\} \cup N(x)$ is a 2-sparse 7-set. So assume $u_1$ is adjacent to at least three vertices in $N(x)$, say without loss of generality $v_1, v_2, v_3$. Note that the sets $\{u_2, u_3, u_4, v_1, v_2, v_3\} \subseteq N(u_1)$ and $N(x)$ are independent since $G$ is triangle-free. Accordingly, we claim that the $\{u_2, u_3, v_1, v_2, v_3, v_4, v_5\}$ shown in Figure \ref{fig:proof_t_2_7} is a 2-sparse 7-set. Indeed, the graph induced by this set can only have edges between vertices in $\{v_4, v_5\}$ and $\{u_2, u_3\}$; yielding at most two neighbors for any vertex. Consequently, there exists a 2-sparse set of size 7 in every triangle-free graph of order 13. Hence, the desired result.

\begin{figure}[htb]
\centering
\begin{tikzpicture}

  \node[draw, circle, inner sep=2pt] (x) at (0,2.5) {};
  \node[anchor=north, xshift=-2pt, yshift=-2pt] at (0,2.5) {$x$};
  
  \foreach \x in {6,5,4,3,2,1}{
   \node[draw, circle, inner sep=2pt] (L\x) at (2,6-\x) {};
   \node[anchor=north, yshift=-2.8pt] at (2,6-\x) {$v_\x$};
   }
   
  \foreach \x in {6,5,4,3,2,1}{
    \node[draw, circle, inner sep=2pt] (R\x) at (4,6-\x) {};
   \node[anchor=north, yshift=-2.8pt] at (4,6-\x) {$u_\x$}; 
    }
  \foreach \i in {3,2,1}{
      \draw (L\i) -- (R1);
  }
   \foreach \i in {6,5,4,3,2,1}{
      \draw (L\i) -- (x);
  }
  \draw (R1) .. controls (4.3,4.8) and (4.3,4.2).. (R2);
  \draw (R1) .. controls (4.75,4.4) and (4.75,3.6).. (R3);
  \draw (R1) .. controls (5.2,4.2) and (5.2,2.8).. (R4);

\draw [yshift=-0.2cm](1.5,5.5) rectangle (2.5,0.5);
\draw [yshift=-0.2cm](3.5,4.5) rectangle (4.5,2.5);

\end{tikzpicture}
\caption{Illustration for the proof of $T_2(7)=13$ in Theorem \ref{thm:T_2(7)=13}.}
\label{fig:proof_t_2_7}
\end{figure}
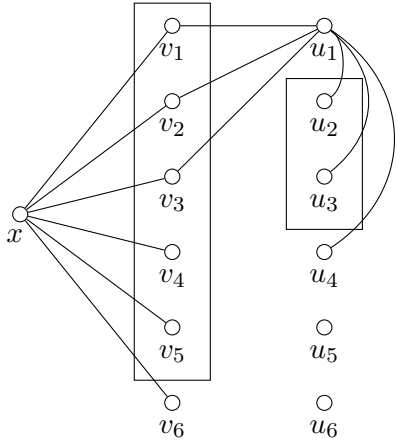
\end{proof}


As for $k=3$, the first non-trivial value $T_3(5)=6$ is implied by Theorem \ref{thm:T_k(k+2)}. The next value for all $k\geq3$, namely $T_k(k+3)$, is provided by the following theorem. For $k=3$, it implies $T_3(6)=8$.

%
%
%
%
%
%

\begin{thm}\label{thm:T_k(k+3)}
With the preceeding notation, $T_k(k+3)=k+5$, for $k \geq 3$ with $K_{2,k+2}$ and $K_{2,k+2}$ with a missing edge as the only two extremal graphs.
\end{thm}
\begin{proof}
$K_{2,k+2}$ and $K_{2,k+2}$ with a missing edge are two graphs which do not include any triangles nor any $k$-sparse $(k+3)$-sets. Therefore, $T_k(k+3) \geq k+5$.

Let $G$ be a triangle-free graph of order $k+5$ and let $x$ be a vertex of maximum degree. If $\Delta(G) \leq k$ then $G$ is $k$-sparse. If $\Delta(G) \geq k+3$, then $N(x)$ is independent and hence contains a $k$-sparse $(k+3)$-set.

Suppose $\Delta(G) = k+2$, and let $y_1$ and $y_2$ be the vertices in $V \setminus N[x]$. If $y_1y_2\notin E$ then, since $k \geq 3$, the set $\{x,y_1,y_2\}$ together with any $k$ vertices from $N(x)$ is a $k$-sparse $(k+3)$-set. If $y_1y_2 \in E$ then no vertex $v \in N(x)$ is adjacent to both $y_1$ and $y_2$, since $G$ is triangle-free. Note, $|N(x)| = k+2$, it follows that at least one of $y_1$ and $y_2$ has at most $k$ neighbors in $N(x)$, say $y_1$. Then $\{y_1\} \cup N(x)$ is a $k$-sparse $(k+3)$-set.

So suppose $\Delta(G) = k+1$, and let $y_1,y_2$ and $y_3$ be the vertices in $V \setminus N[x]$. Since $G$ is triangle-free, there exist a non-edge between two of these three vertices. Without loss of generality, say $y_1y_2 \notin E$. Since $k \geq 3$, the set $\{x, y_1, y_2\}$ together with any $k$ vertices from $N(x)$ is a $k$-sparse $(k+3)$-set. Thus, $T_k(k+3)=k+5$.

Now, let us show that there is no extremal graph other than $K_{2,k+2}$ and $K_{2,k+2}$ with a missing edge. So suppose $G$ is an extremal graph of order $k+4$ and $G$ is not $K_{2,k+2}$ nor $K_{2,k+2}$ minus an edge. Let again $x$ be a vertex of maximum degree. It follows from our previous observations that we can assume $k+1  \leq \Delta(G) \leq k+2$ or else there is a trivial $k$-sparse $(k+3)$-set, contradicting the fact that $G$ is extremal.

So suppose $\Delta(G)=k+1$.  Let $u$ and $w$ be the two vertices not in $N[x]$. Suppose $uw\notin E$, then removing an element of $N(x)$ from $G$ creates a $k$-sparse $(k+3)$-set.  So suppose $uw\in E$.  
Then $u$ and $w$ share no common neighbors since $G$ is triangle-free. Thus one of them, say $u$, is adjacent with at most $k-1$ vertices in $N(x)$.  If $w$ is adjacent with some vertex, say $y$, in $N(x)$ then removing $y$ from $G$ produces a $k$-sparse $(k+3)$-set.  If $w$ is not adjacent with anything in $N(x)$ then it has degree 1.  Thus, removing $x$ produces a $k$-sparse $(k+3)$-set. 

So suppose $\Delta(G)=k+2$.  Let $u$ be the vertex not in $N[x]$.  If the degree of $u$ is $k+1$ or $k+2$ then $G$ is a graph forbidden above.  So suppose the degree of $u$ is at most $k$.  Then removing $x$ produces a $k$-sparse $(k+3)$-set, completing the proof.
\end{proof}
Next, we prove the analogous of Theorems \ref{thm:T_k(k+2)} and \ref{thm:T_k(k+3)} for $k\geq 4$.

\begin{thm}\label{thm:T_k(k+4)}
With the preceeding notation, $T_k(k+4)=k+7$, for $k \geq 4$ with $K_{3,k+3}$ as an extremal graph. 
\end{thm}
\begin{proof}
Observe that $K_{3,k+3}$ does not include any triangles nor any $k$-sparse $(k+4)$-sets. Therefore, $T_k(k+3) \geq k+7$.

Let $G$ be a triangle-free graph of order $k+7$ and let $x$ be a maximum degree vertex of $G$. Similarly to the proof of Theorem \ref{thm:T_k(k+3)}, if $\Delta(G) \geq k+4$ then $N(x)$ includes a $k$-sparse $(k+4)$-set. If $\Delta(G) \leq k$ then $G$ is $k$-sparse.

Suppose $\Delta(G) = k+3$, and denote the remaining vertices in $V \setminus N[x]$ by $Y = \{y_1, y_2, y_3\}$. If $Y$ is an independent set, then $G$ is bipartite with $\{x\} \cup Y$ as one independent set and $N(x)$ as the other. Consequently, $k$ vertices from $N(x)$ together with $\{x\} \cup Y$ is a $k$-sparse $(k+4)$-set, since $k \geq 4$. If $Y$ is not an independent set, say $y_1y_2 \in E$, then a vertex from $N(x)$ cannot be adjacent to both $y_1$ and $y_2$, since $G$ is triangle-free. For $k \geq 4$, either $y_1$ or $y_2$ has at most $k-1$ neighbors in $N(x)$, say it is $y_1$. Then $\{y_1\} \cup N(x)$ is a $k$-sparse $(k+4)$-set.

So suppose $\Delta(G) = k+2$. Denote the vertices in $V \setminus N[x]$ by $Y = \{y_1, y_2, y_3, y_4\}$. If $\alpha(G[Y]) \geq 3$, noting $k \geq 4$, an independent set of $G[Y]$ of size 3, together with $x$ and $k$ vertices from $N(x)$ is a $k$-sparse $(k+4)$-set. If $\alpha(G[Y]) \leq 2$, then since $G$ is triangle-free $G[Y]$ has $2K_2$ as a subgraph. Without loss of generality, say  $\{y_1y_2,y_3y_4\} \in E$. Similar to the reasoning before, a vertex from $N(x)$ cannot be adjacent to both $y_1$ and $y_2$. Consequently, at least one of $y_1$ and $y_2$ has at most $k-1$ neighbors in $N(x)$, say $y_1$ has this property. By symmetry, we can also assume $y_3$ has at most $k-1$ neighbors in $N(x)$. Then, $\{y_1,y_3\} \cup N(x)$ is a $k$-sparse $(k+4)$-set.

So suppose $\Delta(G) = k+1$, and denote the vertices in $V \setminus N[x]$ by $Y = \{y_1, y_2, y_3, y_4, y_5\}$. If $G[Y]$ is bipartite, then $\alpha(G[Y]) \geq 3$, since $k \geq 4$, an independent set of size 3, together with $x$ and $k$ vertices from $N(x)$ is a $k$-sparse $(k+4)$-set. If $G[Y]$ is not bipartite, then it induces a $C_5$ with vertices $y_1, y_2, y_3, y_4, y_5$ in order, since $G$ is triangle-free. Note a vertex of $C_5$ can have at most $k-1$ neigbors in $N(x)$, since $\Delta(G) = k+1$. As a result, for $k \geq 4$, the set $\{y_1,y_2,y_4, x\}$ and $k$ vertices from $N(x)$ is a $k$-sparse $(k+4)$-set.
\end{proof}

We stop proving exact values of $T_k(j)$ and leave the computation of further values using a computer enumeration algorithm for Section \ref{sec:computer}. 
We conclude this section with the following result that establishes the growth rate of $T_k(n)$.
\begin{thm}\label{T_k(n)}
For fixed $k$, we have $T_k(n) = \Theta(\frac{n^2}{\log n})$.
\end{thm}
\begin{proof} 
Fix $k\geq 1$. We know that $c_1 \frac{n^2}{\log n} \leq R(3,n) \leq c_2 \frac{n^2}{\log n}$ for some positive constants $c_1$ and
$c_2$. The first bound is established in \cite{kim} 
and the second in \cite{ajtai}. Note, if $G$ is a triangle-free graph of order at least $ c_2 \frac{n^2}{\log n}$ then $G$ contains an independent set of order $n$. Thus, it contains a $k$-sparse set of order $n$. Accordingly, $T_k(n) \leq c_2 \frac{n^2}{\log n}$ .


So set $j=\Bigl\lfloor c_1 \frac{n^2}{\log n}\Bigr\rfloor$ and let $H$ be a triangle-free graph of order $j$ which contains no independent set of order $n$. Let $H’$ be the lexicographic product of $H$ with an empty graph of order $2k$. Informally, we can think of blowing up each vertex of $H$ with $2k$ isolated vertices, while preserving adjacencies. For a vertex $v$ in $H$, let $S_v$ be the “blown up” vertices of $H’$ that correspond with $v$. Thus, for each $v$, the set $S_v$ is independent and if $uv$ is an edge of $H$, then every vertex of $S_u$ is adjacent with each vertex of $S_v$. Further, an independent set having order $n$ in $H$ corresponds with an independent set in $H’$ with order $2kn$. Note also that $H’$ contains no triangle.

Let $T$ be a $k$-sparse set of $H’$ having maximum order. Note, $|T|\geq 2k\alpha(H)$. Suppose $u$ and $v$ are adjacent in $H$ and $T$ includes vertices from both $S_u$ and $S_v$. Say $x$ and $y$ are in $T$ with $x \in S_u$ and $y \in S_v$. Note, at most $k$ elements of $S_u$
belong to $T$. Further, $x$ can be adjacent to at most $k$ elements of $T$. Now, remove all members of $T$ adjacent to $x$ and replace them with all vertices in $S_u$. The new set is a $k$-sparse set of $H'$ also of maximum order. Repeating this operation for every adjacent pair $x$ and $y$ in $T$ gives an independent set of $H'$ of the same order. Notice, when the vertices of $T$ are ``shrunk" to $H$, an independent set of $H$ is formed. Accordingly, $|T|=\alpha_k(H')=2k\alpha(H)$. As $\alpha(H)<n$, we note that $H’$ is a graph of order $2kj$ which contains no $k$-sparse set of order $2kn$.

Thus, $T_k(2kn)>2kj$. As $k$ is fixed we are allowed a change of variable and note $T_k(n) \geq c_3 \frac{n^2}{\log n}$, for some positive constant $c_3$, and thus our desired conclusion.
\end{proof}

\section{Some defective Ramsey numbers in triangle-free graphs}\label{sec:R}

As noted earlier, the only interesting 1-defective Ramsey numbers in triangle-free graphs are $R_1^{\Delta}(3,j)$ for $j\geq 3$, and $R_1^{\Delta}(4,j)$ for $j\geq 4$ since there is no 1-dense $i$-set for $i\geq 5$ in a triangle-free graph. 

Let us first deal with $R_1^{\Delta}(3,j)$ for $j\geq 3$, and more generally with $R_k^{\Delta}(k+2,j)$ for $j\geq k+2$. It is enough to note that the proof for $R_k(k+2,j)=j$ for all $j\geq k+2$ in general graphs given in \cite{ekim2013some} is also valid in triangle-free graphs. Thus, we have the following, which is also certified by computer enumeration in Section \ref{sec:computer} (see Tables \ref{table:R_{1}(i,j)}, \ref{table:R_{2}(i,j)}, \ref{table:R_{3}(i,j)}, \ref{table:R_{4}(i,j)}).

\begin{rem}
With the preceeding notation, $R_k^{\Delta}(k+2,j)=j$ for $j\geq k+2$.
\end{rem}

In what follows, we investigate $R_1^{\Delta}(4,j)$ for $j\geq 4$.
Recall that $C_4$ is the only triangle-free 1-dense 4-set. In this section, we show $R_1^{\Delta}(4,4)=6, R_1^{\Delta}(4,5)=8, R_1^{\Delta}(4,6)=10$ and $R_1^{\Delta}(4,7)=13$. We provide extremal graphs for each result. Uniqueness will be established in Section \ref{sec:computer} using computer enumeration.
 
\begin{thm}\label{thm:??}
With the preceeding notation, $R_1^{\Delta}(4,4)=6$ with the unique extremal graph being $C_5$.
\end{thm}
\begin{proof}
Note that $C_5$ is triangle-free which does not contain any 1-dense 4-set nor 1-sparse 4-set. Thus, $R_1^{\Delta}(4,4)\geq 6$. Consider a triangle-free graph $G$ of order 6. If $G$ has a $C_4$, then it is a 1-dense 4-set. 
So assume that  $G$ does not contain $C_4$. If $G$ contains $C_5$, the vertex $x$ that is not on the $C_5$ can only be adjacent to a single vertex from the $C_5$, otherwise there would be a triangle or a $C_4$. Denote the vertices on the cycle by $v_1, v_2, v_3, v_4, v_5$ in order, with $x$ being possibly adjacent to one vertex, say wothout loss of generality $v_1$, and no other vertex. Note the set $\{x, v_1, v_3, v_4\}$ is a 1-sparse 4-set, whether or not $x$ is adjacent to $v_1$. So suppose that $G$ has a $C_6$. Then G is a $C_6$ and two opposing edges create a 1-sparse 4 set. Lastly, consider the case where $G$ has no cycles, that is $G$ is a forest. In this case, $G$ is actually a bipartite graph. If it is an unbalanced bipartite graph, meaning that one of the independent sets is of size at least 4, then that set is a 1-sparse 4-set. If $G$ is a balanced bipartite graph, there is a vertex $v$ that has at most 1 neighbor in the other independent set, say $U$, since all forests contain a vertex of degree at most 1.  Then, $\{v\} \cup U$ is a 1-sparse 4-set. In conclusion, every triangle-free graph of order $6$ has either a 1-dense 4-set or a 1-sparse 4-set.  
\end{proof}

\begin{thm}
With the preceeding notation, $R_1^{\Delta}(4,5)=8$ with the unique extremal graph being $C_7$.
\end{thm}
\begin{proof}
Note that $C_7$ is a triangle-free graph of order 7 which does not contain any 1-dense 4-set nor 1-sparse 5-set. Thus, $R_1^{\Delta}(4,5)\geq 8$. Let $G$ be a triangle-free graph of order 8. If it has a 4-cycle, then it has a 1-dense 4-set. If $G$ has no $C_4$, then we will show that it contains a 1-sparse 5-set. 

If $G$ has girth 5, then, the three vertices outside a 5-cycle $C$, denote by $v_1, v_2, v_3$, each can be adjacent to at most 1 vertex from the cycle. Otherwise, a triangle or a $C_4$ would exist. Call $x_1, x_2, x_3, x_4, x_5$ the vertices of $C$ in order. If a vertex from $C$, say  $x_1$, is adjacent to 2 or 3 vertices in $\{v_1, v_2, v_3\}$, then the set $(N(x_1)\setminus C) \cup  \{x_2, x_3, x_5\}$ contains a 1-sparse 5-set. Now, assume that all vertices from $C$ have at most 1 neighbor from $\{v_1, v_2, v_3\}$. Under these conditions, we can choose two vertices in $\{v_1, v_2, v_3\}$, say without loss of generality $v_1$ and $v_2$, such that $N(\{v_1, v_2\})$ do not contain two vertices of $C$ which are adjacent. Then, we can choose a set $A$ in $ C \setminus N(\{v_1, v_2\})$ which is a 1-sparse 3-set. Then, $A \cup \{v_1, v_2\}$ is a 1-sparse 5-set.

Now, assume that $G$ has girth 6 and let $C$ be a 6-cycle. Then, call $v_1$ and $v_2$ the two vertices that are not on $C$. Since the girth is 6, each one of $v_1$ and $v_2$ has at most one neighbor in $C$. So, there exists a set $A$ in $C \setminus N(\{v_1, v_2\})$ that is a 1-sparse 3-set. Then, $A \cup \{v_1, v_2\}$ is a 1-sparse 5-set. Hence, the girth is at least 7.

If $G$ has girth 7, then it is a $C_7$ with vertices $x_1, x_2, x_3, x_4, x_5, x_6, x_7$ and a remaining vertex outside the cycle, say $v$. Indeed, $v$ can be adjacent to at most 1 vertex, say $x_1$, from the $C_7$, since girth is 7. The set $\{v, x_2, x_3, x_5, x_6\}$ is a 1-sparse 5-set whether $v$ is adjacent to $x_1$ or not. Lastly, if $G$ has girth 8 and is a $C_8$ with vertices $x_1, x_2, x_3, x_4, x_5, x_6, x_7, x_8$, then the set $\{x_1, x_3, x_4, x_6, x_7\}$ is a 1-sparse 5 set.

Finally the case where $G$ is a forest remains. If $G$ is a forest, then it is bipartite. If it is an unbalanced bipartite graph, then the independent set with higher size includes a 1-sparse 5-set. Assume that $G$ is a balanced bipartite graph. There must exist a vertex with degree at most 1, since $G$ is also a forest. This pendent vertex and an independent set it does not belong to together create a 1-sparse 5-set. Therefore, every traingle-free graph of order 8 includes either a 1-dense 4-set or a 1-sparse 5-set. 
\end{proof}

\begin{thm}
With the preceeding notation, $R_1^{\Delta}(4,6)=10$ with the unique extremal graph being the graph given in Figure \ref{fig:extremal_r_1_4_6}.
\end{thm}
\begin{proof}
Consider the graph in Figure \ref{fig:extremal_r_1_4_6}. It is a triangle-free graph which does not contain any 1-dense 4-set nor 1-sparse 6-set. So, $R_1^{\Delta}(4,6)\geq 10$. Let $G$ be a triangle-free graph of order 10. If it has a 4-cycle, then it has a 1-dense 4-set. So assume $G$ has no $C_4$, then we will show that it contains a 1-sparse 6-set. 

If $\Delta(G) \geq 6$, then the neighborhood of a maximum degree vertex includes a 1-sparse 6-set. So, assume $\Delta(G)\leq 5$. Suppose $G$ has a vertex $x$ of degree 5. Note every vertex outside $N[x]$ is adjacent to at most 1 vertex from $N(x)$, or else a 4-cycle is formed. Then $N(x)$ and a vertex outside $N_G[x]$ is a 1-sparse 6-set. So assume $\Delta(G)\leq 4$.

Suppose $G$ has a vertex $x$ of degree 4 and let $N(x) = \{v_1, v_2, v_3, v_4\}$ and $U = V \setminus N[x] =  \{u_1, u_2, u_3, u_4, u_5\}$. If $U$ is a 1-sparse 5-set, then $U \cup \{x\}$ is a 1-sparse 6-set. So assume $U$ is not 1-sparse, thus there is a vertex from $U$, say $u_1$, which is adjacent two other vertices in $U$, say $u_2$ and $u_3$, without loss of generality. Note any vertex from $U$ can be adjacent to at most one vertex in $N(x)$, or else a $C_4$ is formed. Moreover, $u_2u_3\notin E$ since $G$ is triangle-free. Likewise, a vertex from $N(x)$ cannot be adjacent to both $u_2$ and $u_3$, otherwise a $C_4$ is induced by that vertex and $\{u_1, u_2, u_3\}$. Consequently, $N(x) \cup \{u_2, u_3\}$ is a 1-sparse 6-set. So, we may assume $\Delta(G)\leq 3$.

If there is a vertex $x$ of degree 1 in $G$, then $V\setminus N[x]$ has a 1-sparse 5-set $S$ by $R_1^{\Delta}(4,5)=8$; thus $\{x\} \cup S$ is a 1-sparse 6-set. So assume every vertex in $G$ has degree at least 2.

Assume there is a vertex $x$ of degree 3 and $U = V \setminus N[x] =  \{u_1, u_2, u_3, u_4, u_5, u_6\}$. If $U$ is 1-sparse, then it is a 1-sparse 6-set. Suppose $U$ is not 1-sparse, thus there exists a vertex in $U$, say $u_1$, such that $d_U(u_1)\geq 2$. If $d_U(u_1)=3$, then $N(x) \cup N(u_1)$ is a 1-sparse 6-set. So assume every vertex $u\in U$ has $d_U(u)\leq 2$. Thus, $U$ induces a collection of disjoint paths and cycles. Note $G$ has no triangle, nor $C_4$. Also, $G[U]$ cannot have a $C_5$ since the remaining vertex in $U$ would be adjacent to two vertices in $N(x)$ (since $\delta(G)\geq 2$) forming a $C_4$. So the only cycle in $G[U]$ can be a 6-cycle. In this case, take an independent set $I$ of 3 vertices in this 6-cycle; then $I\cup N(x)$ is a 1-sparse 6-set (any vertex in $I$ has at most one neighbor in $N(x)$ and vice versa, or else a $C_4$ is formed). So assume $G[U]$ is a collection of paths. Observe every pendant vertex in $G[U]$ has to be adjacent to at least one vertex in $N(x)$ (since $\delta(G)\geq 2$); indeed it is adjacent to exactly one vertex in $N(x)$ or else a $C_4$ is formed. Thus, there is at most 3 pendant vertices in $G[U]$. By the Handshaking Lemma, $G[U]$ has exactly 2 pendant vertices; thus $G[U]$ is a $P_6$. Note there is a vertex in $N(x)$, say $v$, which is not adjacent to the end-vertices of the $P_6$. Clearly, $v$ has at most two neighbors in $P_6$. Moreover, it is possible to chose a 1-sparse 4-set in $P_6\setminus N(v)$ which together with $x$ and $v$ forms a 1-sparse 6-set. So we may assume $G$ is 2-regular. Thus it is either a $C_{10}$ or two disjoint copies of $C_5$. In both cases, there is a 1-sparse 6-set.

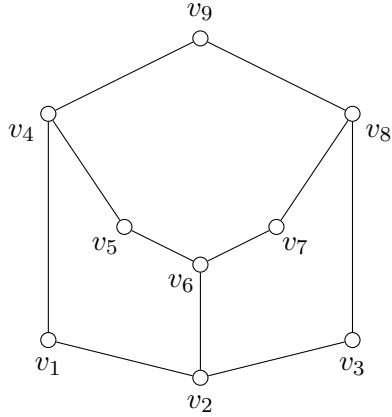
\begin{figure}
\centering
\begin{tikzpicture}

\foreach \x/\y/\name/\xsh/\ysh in {
0/0/1/0/-0.1,
2/-0.5/2/0/-0.1,
4/0/3/0/-0.1,
0/3/4/-0.35/0,
1/1.5/5/-0.25/0,
2/1/6/-0.25/0,
3/1.5/7/0.25/0,
4/3/8/0.35/0,
2/4/9/0/0.6
} {
  \node[draw, circle, inner sep=2pt] (v\name) at (\x, \y) {};
  \node[anchor=north, xshift=\xsh cm, yshift=\ysh cm] at (\x, \y) {$v_\name$};
}

\foreach \from/\to in {1/2, 2/3, 3/8, 4/5, 5/6, 6/7, 7/8, 8/9, 1/4, 4/9, 2/6} {
  \draw (v\from) -- (v\to);
}

\end{tikzpicture}
\caption{The unique extremal graph for $R_1^\Delta(4,6)=10$.}
\label{fig:extremal_r_1_4_6}
\end{figure}

\end{proof}

\begin{thm}
With the preceeding notation, $R_1^{\Delta}(4,7)=13$ with exactly 2 extremal graphs given in Figure \ref{fig:extremal_r_1_4_7}.
\end{thm}
\begin{proof}
The graphs in Figure \ref{fig:extremal_r_1_4_7} are both triangle-free graphs which do not contain any 1-dense 4-set nor 1-sparse 7-set. So, $R_1^{\Delta}(4,7)\geq 13$. 

Suppose to the contrary, there is some triangle-free graph of order 13 containing no 1-dense 4-set nor 1-sparse 7-set. Let $G$ be such a graph. Note $G$ contains no 4-cycle induced or otherwise.

If $\delta(G)\leq 2$ then remove a minimum degree vertex $x$ along with its neighborhood. Note a graph on at least 10 vertices remains. By $R_1^{\Delta}(4,6)=10$, the remaining graph has a 1-dense 4-set, a contradiction; or a 1-sparse 6-set, which together with $x$, forms a 1-sparse 7-set. So, assume every vertex has degree at least 3.

If $\Delta(G) \geq 7$, then the neighborhood of a maximum degree vertex includes a 1-sparse 7-set. So we may assume $\Delta(G) \leq 6$. If $\Delta(G) = 6$, let $x$ be a vertex of degree 6 and consider a vertex $y\in V\setminus N[x]$. If $y$ has two neighbors in $N(x)$, then these two neighbors together with $x$ and $y$ form a $C_4$. So assume $y$ has at most one neighbor in $N(x)$, then $N(x)\cup \{y\}$ is a 1-sparse 7-set. So we may assume $\Delta(G) \leq 5$.

So suppose $G$ has a vertex $x$ of degree 5. Then $V\setminus N[x]$ induces  a triangle-free graph on 7 vertices; since $R(3,3)=6$, it has an independent set $A$ of size 3. If no vertex of $A$ is adjacent to more than one vertex in $N(x)$, then $A\cup N(x)$ contains a 1-sparse 7-set. Note no vertex of $A$ is adjacent to more than one vertex in $N(x)$, or else $C_4$ is present. If all three vertices of $A$ are adjacent to the same vertex $y\in N(x)$, then $A\cup (N(x)\setminus y)$ is a 1-sparse 7-set. If there are two vertices, say $u,v\in A$ are adjacent to the same vertex of $N(x)$. Then $(A\setminus \{u\}\cup N(x))$ is a 1-sparse 7-set. If every vertex in $N(x)$ is adjacent to at most one vertex in $A$, $A\cup N(x)$ contains a 1-sparse 7-set. Thus, $G$ has no vertex of degree 5.

So, suppose the maximum degree of $G$ is 4.  Let $x$ be a vertex of degree four.  Suppose also that $G$ has a second vertex, say $y$, of degree 4.  Let us consider the case where $x$ and $y$ are non-adjacent.  Note, $x$ and $y$ cannot have two common neighbors, for otherwise $G$ contains a 4-cycle.  So, $x$ and $y$ have at most one common neighbor and hence, $N(x)\cup N(y)$ is a 1-sparse set on at least 7 vertices.  Thus, all vertices of degree 4 are adjacent with $x$.

Pick $y$, a non-neighbor of $x$.  Note, $y$ has degree  exactly 3 since $\delta(G)\geq 3$.  Suppose $N(x)$ and $N(y)$ don’t meet.  Then their union is a 1-sparse 7-set because of the absence of 4-cycles.  So let us assume these neighborhoods meet and $z$ belongs to both.  We note there can be no other vertex belonging to both.  As $G$ contains no vertices of degree 2, we note $z$ is adjacent to some other vertex and this vertex is outside $N[x] \cup N[y]$.  Call one such vertex $w$.  Note, $w$ cannot be adjacent with anything in $N(x)\cup N(y)$ other than $z$, for otherwise a 4-cycle is present in $G$.  Thus, $N(x) \cup N(y) \cup \{w\}$ is a 1-sparse 7-set.  Thus, $G$ contains no vertex of degree 4.

Accordingly, $G$ is 3-regular.  But this is impossible; by the Handshaking Lemma, there is no 3-regular graph of order 13. 

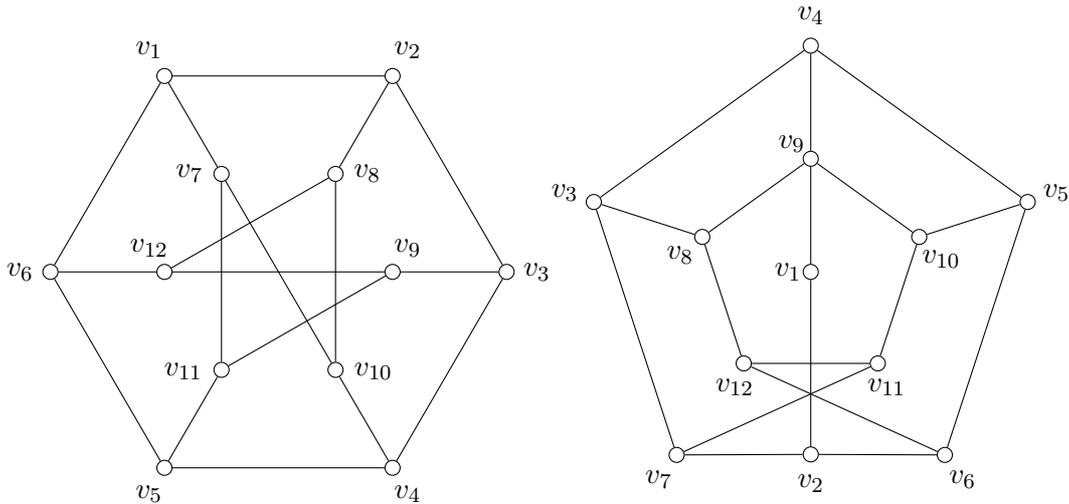
\begin{figure}[htb]
\centering
\begin{tikzpicture}

\foreach \name/\polar/\dist/\xsh/\ysh/\dsh in {
1/120/3.0/0.0/0.0/0.4,
2/60/3.0/0.0/0.0/0.4,
3/0/3.0/0.0/0.0/0.4,
4/-60/3.0/0.0/0.0/0.4,
5/-120/3.0/0.0/0.0/0.4,
6/180/3.0/0.0/0.0/0.4,
7/120/1.5/-0.32/-0.16/0.2,
8/60/1.5/0.32/-0.16/0.2,
9/0/1.5/0.0/0.32/0.2,
10/-60/1.5/0.4/0.16/0.2,
11/-120/1.5/-0.4/0.16/0.2,
12/180/1.5/0.0/0.32/0.2
} {
  \node[draw, circle, inner sep=2pt] (v\name) at (\polar: \dist) {}; 
  \node[xshift=\xsh cm, yshift=\ysh cm] at (\polar: \dist + \dsh) {$v_{\name}$};
}


\foreach \from/\to in {1/2, 2/3, 3/4, 4/5, 5/6, 1/6, 1/7, 2/8, 3/9, 4/10, 5/11, 6/12, 12/9, 9/11, 11/7, 7/10, 10/8, 8/12} {
  \draw (v\from) -- (v\to);
}

\foreach \name/\polar/\dist/\xsh/\ysh/\dsh in {
1/90/0/-0.3/0.0/0,
2/270/2.42/0/0/0.38,
3/162/3.0/0.0/0.0/0.4,
4/90/3.0/0.0/0.0/0.4,
5/18/3.0/0.0/0.0/0.4,
6/306/3.0/0.0/0.0/0.4,
7/234/3.0/0.0/0.0/0.4,
8/162/1.5/-0.1/-0.3/0.2,
9/90/1.5/-0.25/0.0/0.2,
10/18/1.5/0.1/-0.35/0.2,
11/306/1.5/0.0/-0.15/0.2,
12/234/1.5/0.0/-0.15/0.2
} {
  \node[draw, circle, inner sep=2pt, xshift = 7 cm] (v\name) at (\polar: \dist) {}; 
  \node[xshift={\xsh cm + 7 cm}, yshift=\ysh cm] at (\polar: \dist + \dsh) {$v_{\name}$};
}


\foreach \from/\to in {1/2, 1/9, 2/7, 2/6, 3/4, 3/7, 3/8, 4/5, 4/9, 5/10, 5/6, 11/12, 8/9, 9/10, 10/11, 8/12, 7/11, 6/12} {
  \draw (v\from) -- (v\to);
}

\end{tikzpicture}
\caption{The two extremal graphs for $R_1^\Delta(4,7)=13$.}
\label{fig:extremal_r_1_4_7}
\end{figure}

\end{proof}

\section{Computer enumeration} \label{sec:computer}

In this work, we obtain several defective Ramsey numbers in triangle-free graphs using proofs ``by hand". Whenever classical proof techniques hit limits due to the highly combinatorial nature of the extremal graphs and Ramsey numbers, we also make use of a computer based search. We use an adaptation of the Algorithm Sub-extremal given in \cite{SubExt} for triangle-free graphs as described here in Algorithm \ref{algo:SubExtTfree}. Algorithm \ref{algo:SubExtTfree} computes new defective Ramsey numbers and enumerates related extremal graphs for triangle-free graphs. It also serves as a checking mechanism for the proofs made by hand in earlier sections. All of our codes and the extremal graphs we obtain are available online at \cite{github}.

Let us denote by $\mathcal{T}_n^\Delta(k,i,j)$ the set of all triangle-free graphs of order $n$ containing no $k$-dense $i$-set nor $k$-sparse $j$-set. Given $\mathcal{T}_n^\Delta(k,i,j)$, we call a $k$-dense $i$-set or a $k$-sparse $j$-set \textit{a forbidden $k$-defective set}. Note that the set of all extremal graphs for ${R}_k^\Delta(i,j)$ is the set $\mathcal T_{n}^\Delta(k,i,j)$ for $n={R}_k^\Delta(i,j)-1$. Accordingly, a graph in  $\mathcal{T}_{n}^\Delta(k,i,j)$ for $n < {R}_k^\Delta(i,j)-1$ is  called a \textit{sub-extremal} graph for ${R}_k^\Delta(i,j)$.

Algorithm Sub-extremal in \cite{SubExt} computes the defective Ramsey number ${R}_k^{\mathcal G}(i,j)$ and all its extremal graphs for some parameters $k,i,j$ such that $i,j \geq k+2$ and for some graph class $\mathcal G$. In its generic form, it checks whether the generated graphs belong to the desired graph class $\mathcal G$ or not at the very end of the algorithm, and eliminates those not in $\mathcal G$. Our adaptation given in Algorithm \ref{algo:SubExtTfree} differs from Algorithm Sub-extremal only in checking the presence of triangles. Rather than checking for this at the end, we ensure that the generated graphs are always triangle-free by carefully adding each new vertex in Line \ref{line:trianglefree}. 

Algorithm \ref{algo:SubExtTfree} is based on the fact that being triangle-free, ($k$-dense $i$-set)-free, and ($k$-sparse $j$-set)-free are hereditary properties. Given a sub-extremal graph $G$ of order $n$, all graphs of order $n+1$ that have $G$ as an induced subgraph are produced by adding a new vertex with all possible adjacency combinations to the vertices of $G$. If a new graph created by this procedure is also triangle-free and contains no forbidden $k$-defective set for $\mathcal{T}_{n+1}^\Delta(k,i,j)$, it belongs to the set of (sub-)extremal graphs of order $n+1$. Taking $\mathcal{T}_n^\Delta(k,i,j)$ as input, the set $\mathcal{T}_{n+1}^\Delta(k,i,j)$ is generated by this method. We start with the one vertex graph $K_1$ as input. We run Algorithm \ref{algo:SubExtTfree} iteratively giving the output of one iteration as the input of the next iteration. We stop when the output set is empty and declare ${R}_k^\Delta(i,j) = n$ where $n$ is the number of vertices for which the algorithm does not return a graph. This means all triangle-free graphs with the present order $n$ (or larger) contain either a $k$-dense $i$-set or a $k$-sparse $j$-set. We conclude that the last non-empty output set of graphs with ${R}_k^\Delta(i,j)-1$ vertices is the complete list of all extremal graphs for ${R}_k^\Delta(i,j)$.

	\begin{algorithm}[htb]
			\DontPrintSemicolon
		
		\KwInput{$\mathcal{T}_n^\Delta(k,i,j)$ for some $k,i,j$ such that $i,j \geq k+2$}
		\KwOutput{$\mathcal{T}_{n+1}^\Delta(k,i,j)$}
		
		Let $K=\emptyset$. \\ 
		\ForEach {\label{line:input_graph_loop} $G\in \mathcal{T}_n^\Delta(k,i,j)$}{
			\ForAll {$S\subseteq V(G)$}{ \label{line:S}
				\If {$S$ is an independent set} {	\label{line:trianglefree} Take the graph $G_S$ that is formed by adding a new vertex $v$ into $G$ that is adjacent to all vertices in $S$.\\
				Let $add$ = \textbf{TRUE}\\
				\ForAll {\label{line:def_start}$I\subseteq V(G_S)$ such that $v\in I$ and $|I|\in\{i,j\}$}  { 
					\If{\label{line:dense_start}$|I|=i$ and $G[I]$ is $k$-dense}
						{$add$ = \textbf{FALSE} and \textbf{BREAK}}
					\label{line:dense_end}
					\If{$|I|=j$ and $G[I]$ is $k$-sparse}{$add$ = \textbf{FALSE} and \textbf{BREAK} \label{line:def_end}}
				} 
				\If{$add$ = \textbf{\textup{TRUE}}}{
				 	Add $G_S$ into $K$.
				}
				}
		}}
		Return a maximal non-isomorphic set of graphs in $K$.
		\caption{Sub-extremal for Triangle-free Graphs} \label{algo:SubExtTfree}
	\end{algorithm}

Let $G$ be a triangle-free graph of order $n$. Adding a new vertex that is adjacent to an independent set of $G$ creates a new triangle-free graph of order $n + 1$ which contains $G$. By applying this procedure for every independent set of $G$, we obtain all triangle-free graphs of order $n + 1$ containing $G$. This procedure is executed in line \ref{line:trianglefree} of Algorithm \ref{algo:SubExtTfree}. The significant advantage of this method is that it ensures that the new graph is triangle-free without requiring an explicit check for the existence of triangles.



Having guaranteed the absence of triangles, all we need to check is whether one of the forbidden $k$-defective sets is formed. Indeed, since the input graphs have no forbidden $k$-defective sets, if a newly generated graph $G_S$ contains a forbidden $k$-defective set, then this must contain the new vertex $v$. Accordingly, it is sufficient to check all subsets including the new vertex $v$ for forbidden $k$-defective sets in lines \ref{line:def_start} to \ref{line:def_end}. Checking the existence of a $k$-dense $i$-set, in lines \ref{line:dense_start} to \ref{line:dense_end}, is included in the search for $R_k^{\Delta}(4,j)$ values. However, the $k$-dense set checking mechanism is omitted for $T_k(j)$ values which only consider sparse sets.

The nature of Algorithm \ref{algo:SubExtTfree} allows for parallel computing. Thus we implemented both the graph generation and isomorphism checks to execute in parallel and obtained improved runtime efficiency. In the graph generation, each thread works with a separate graph from the set of input graphs in line \ref{line:input_graph_loop}. All generated and valid graphs are pooled together in an array. In this pool, isomorphic copies of graphs exist and removing isomorphic copies of a graph is a challenge for this algorithm. The isomorphism checks are carried out by comparing graphs using their canonical labelings which are calculated using the \emph{nauty} program \cite{mckaynauty}. This isomorphism check is also programmed to execute in parallel with each thread checking a different graph and utilizing mutex locks which allows the threads to safely access and modify the shared data structures. 

The program was implemented in C++ and executed on a personal computer with 8 gigabytes of RAM and Apple M1 chip which has 8 cores and maximum CPU clock rate of 3.2 GHz. For the relatively small defective Ramsey numbers, the runtime is trivially quick. However, as the graph size increases, runtimes grow exponentially both in generation and isomorphism check phases. The longest runtime encountered for a defective Ramsey number, which is $R_4^\Delta(9,11) = 18$, was approximately 2.5 hours.

We found several defective Ramsey numbers that we cannot prove by hand as well as the number of the extremal graphs using Algorithm \ref{algo:SubExtTfree}. Table \ref{table:T_k(j)} displays $T_k(j)$ values computed by Algorithm \ref{algo:SubExtTfree} and the corresponding number of extremal graphs for each number computed. Tables \ref{table:R_{1}(i,j)}, \ref{table:R_{2}(i,j)}, \ref{table:R_{3}(i,j)}, \ref{table:R_{4}(i,j)} display similar results obtained for defective Ramsey numbers $R_{k}^\Delta(i,j)$ for $k = 1,2,3$ and $4$. Missing numbers in these tables could not be obtained due to insufficient memory necessary for storing the sub-extremal graphs. The subsequent data offers an insight into the memory needs of the defective Ramsey numbers we computed and those we could not determine. In the computation of $T_1(7)=18$, the maximum number of sub-extremal graphs encountered is 1243785 (for order 13) and this number is then reduced to 108243 non-isomorphic graphs. In the computation of $T_1(8)$, 8958224 sub-extremal graphs of order 12 are reduced to 822971 non-isomorphic graphs. In the generation of sub-extremal graphs of order 13, the program was terminated by the operating system. Optimizing the implementation of Algorithm \ref{algo:SubExtTfree} and utilizing better hardware may result in the computation of a few more defective Ramsey numbers. However, due to the exponential growth in the number of graphs we need to construct in Line~\ref{line:S} of Algorithm~\ref{algo:SubExtTfree}, increase in the runtimes and memory requirements will remain as two significant challenges. 

\begin{table}[!htb]
\centering
\caption{$T_k(j)$ values and corresponding number of extremal graphs.}
\label{table:T_k(j)}\small
\resizebox{\textwidth}{!}{
\begin{tabular}{ll|llllllllll|}
\cline{3-12}
                                            &   & \multicolumn{10}{c|}{$j$}                                                                                                                                                                                                                                                                                                                                                                                                                                                                                  \\ \cline{3-12} 
                                            &   & \multicolumn{1}{l|}{3}     & \multicolumn{1}{l|}{4}                            & \multicolumn{1}{l|}{5}                             & \multicolumn{1}{l|}{6}                             & \multicolumn{1}{l|}{7}                             & \multicolumn{1}{l|}{8}                              & \multicolumn{1}{l|}{9}                              & \multicolumn{1}{l|}{10}                              & \multicolumn{1}{l|}{11}                              & 12                               \\ \hline
\multicolumn{1}{|l|}{}                      & 1 & \multicolumn{1}{l|}{5 (1)} & \multicolumn{1}{l|}{7 (2)}                        & \multicolumn{1}{l|}{11 (1)}                        & \multicolumn{1}{l|}{13 (16)}                       & \multicolumn{1}{l|}{18 (1)}                        & \multicolumn{1}{l|}{}                               & \multicolumn{1}{l|}{}                               & \multicolumn{1}{l|}{}                                & \multicolumn{1}{l|}{}                                &                                  \\ \cline{2-12} 
\multicolumn{1}{|l|}{}                      & 2 & \multicolumn{1}{l|}{3 (2)} & \multicolumn{1}{l|}{\cellcolor[HTML]{FE996B}5(1)} & \multicolumn{1}{l|}{9 (2)}                         & \multicolumn{1}{l|}{11 (6)}                        & \multicolumn{1}{l|}{13 (288)}                      & \multicolumn{1}{l|}{16 (281)}                       & \multicolumn{1}{l|}{}                               & \multicolumn{1}{l|}{}                                & \multicolumn{1}{l|}{}                                &                                  \\ \cline{2-12} 
\multicolumn{1}{|l|}{}                      & 3 & \multicolumn{1}{l|}{3 (2)} & \multicolumn{1}{l|}{4 (3)}                        & \multicolumn{1}{l|}{\cellcolor[HTML]{FE996B}6 (1)} & \multicolumn{1}{l|}{\cellcolor[HTML]{F8FF00}8 (2)} & \multicolumn{1}{l|}{13 (5)}                        & \multicolumn{1}{l|}{15 (40)}                        & \multicolumn{1}{l|}{17 (9713)}                      & \multicolumn{1}{l|}{}                                & \multicolumn{1}{l|}{}                                &                                  \\ \cline{2-12} 
\multicolumn{1}{|l|}{}                      & 4 & \multicolumn{1}{l|}{3 (2)} & \multicolumn{1}{l|}{4 (3)}                        & \multicolumn{1}{l|}{5 (7)}                         & \multicolumn{1}{l|}{\cellcolor[HTML]{FE996B}7 (1)} & \multicolumn{1}{l|}{\cellcolor[HTML]{F8FF00}9 (2)} & \multicolumn{1}{l|}{\cellcolor[HTML]{34CDF9}11 (7)} & \multicolumn{1}{l|}{17 (19)}                        & \multicolumn{1}{l|}{19 (606)}                        & \multicolumn{1}{l|}{}                                &                                  \\ \cline{2-12} 
\multicolumn{1}{|l|}{}                      & 5 & \multicolumn{1}{l|}{3 (2)} & \multicolumn{1}{l|}{4 (3)}                        & \multicolumn{1}{l|}{5 (7)}                         & \multicolumn{1}{l|}{6 (14)}                        & \multicolumn{1}{l|}{\cellcolor[HTML]{FE996B}8 (1)} & \multicolumn{1}{l|}{\cellcolor[HTML]{F8FF00}10 (2)} & \multicolumn{1}{l|}{\cellcolor[HTML]{34CDF9}12 (7)} & \multicolumn{1}{l|}{\cellcolor[HTML]{32CB00}14 (46)} & \multicolumn{1}{l|}{21 (112)}                        &                                  \\ \cline{2-12} 
\multicolumn{1}{|l|}{}                      & 6 & \multicolumn{1}{l|}{3 (2)} & \multicolumn{1}{l|}{4 (3)}                        & \multicolumn{1}{l|}{5 (7)}                         & \multicolumn{1}{l|}{6 (14)}                        & \multicolumn{1}{l|}{7 (38)}                        & \multicolumn{1}{l|}{\cellcolor[HTML]{FE996B}9 (1)}  & \multicolumn{1}{l|}{\cellcolor[HTML]{F8FF00}11(2)}  & \multicolumn{1}{l|}{\cellcolor[HTML]{34CDF9}13 (7)}  & \multicolumn{1}{l|}{\cellcolor[HTML]{32CB00}15 (46)} & \cellcolor[HTML]{FFCCC9}17 (723) \\ \cline{2-12} 
\multicolumn{1}{|l|}{\multirow{-7}{*}{$k$}} & 7 & \multicolumn{1}{l|}{3 (2)} & \multicolumn{1}{l|}{4 (3)}                        & \multicolumn{1}{l|}{5 (7)}                         & \multicolumn{1}{l|}{6 (14)}                        & \multicolumn{1}{l|}{7 (38)}                        & \multicolumn{1}{l|}{8 (107)}                        & \multicolumn{1}{l|}{\cellcolor[HTML]{FE996B}10 (1)} & \multicolumn{1}{l|}{\cellcolor[HTML]{F8FF00}12 (2)}  & \multicolumn{1}{l|}{\cellcolor[HTML]{34CDF9}14 (7)}  & \cellcolor[HTML]{32CB00}16 (46)  \\ \hline
\end{tabular}
}
\end{table}

\begin{table}[!htb]
\centering
\caption{$R_{1}^{\Delta}(i,j)$  values and corresponding number of extremal graphs.}
\label{table:R_{1}(i,j)}
\begin{tabular}{ll|lllllll|}
\cline{3-9}
                                           &   & \multicolumn{7}{c|}{$j$}                                                                                                                                                                \\ \cline{3-9} 
                                           &   & \multicolumn{1}{l|}{3}     & \multicolumn{1}{l|}{4}     & \multicolumn{1}{l|}{5}     & \multicolumn{1}{l|}{6}      & \multicolumn{1}{l|}{7}      & \multicolumn{1}{l|}{8}      & 9      \\ \hline
\multicolumn{1}{|l|}{\multirow{2}{*}{$i$}} & 3 & \multicolumn{1}{l|}{3 (2)} & \multicolumn{1}{l|}{4 (2)} & \multicolumn{1}{l|}{5 (3)} & \multicolumn{1}{l|}{6 (3)}  & \multicolumn{1}{l|}{7 (4)}  & \multicolumn{1}{l|}{8 (4)}  & 9 (5)  \\ \cline{2-9} 
\multicolumn{1}{|l|}{}                     & 4 & \multicolumn{1}{l|}{4 (1)} & \multicolumn{1}{l|}{6 (1)} & \multicolumn{1}{l|}{8 (1)} & \multicolumn{1}{l|}{10 (1)} & \multicolumn{1}{l|}{13 (2)} & \multicolumn{1}{l|}{15 (3)} & 18 (4) \\ \hline
\end{tabular}
\end{table}

\begin{table}[!htb]
\centering
\caption{$R_{2}^{\Delta}(i,j)$ values and corresponding number of extremal graphs.}
\label{table:R_{2}(i,j)}
\begin{tabular}{ll|lllllll|}
\cline{3-9}
                                           &   & \multicolumn{7}{c|}{$j$}                                                                                                                                                                 \\ \cline{3-9} 
                                           &   & \multicolumn{1}{l|}{4}     & \multicolumn{1}{l|}{5}     & \multicolumn{1}{l|}{6}     & \multicolumn{1}{l|}{7}      & \multicolumn{1}{l|}{8}       & \multicolumn{1}{l|}{9}      & 10     \\ \hline
\multicolumn{1}{|l|}{\multirow{3}{*}{$i$}} & 4 & \multicolumn{1}{l|}{4 (3)} & \multicolumn{1}{l|}{5 (3)} & \multicolumn{1}{l|}{6 (3)} & \multicolumn{1}{l|}{7 (3)}  & \multicolumn{1}{l|}{8 (3)}   & \multicolumn{1}{l|}{9 (3)}  & 10 (3) \\ \cline{2-9} 
\multicolumn{1}{|l|}{}                     & 5 & \multicolumn{1}{l|}{5 (1)} & \multicolumn{1}{l|}{6 (4)} & \multicolumn{1}{l|}{8 (1)} & \multicolumn{1}{l|}{10 (2)} & \multicolumn{1}{l|}{11 (62)} & \multicolumn{1}{l|}{15 (2)} & 17 (4) \\ \cline{2-9} 
\multicolumn{1}{|l|}{}                     & 6 & \multicolumn{1}{l|}{5 (1)} & \multicolumn{1}{l|}{7 (3)} & \multicolumn{1}{l|}{9 (6)} & \multicolumn{1}{l|}{12 (5)} & \multicolumn{1}{l|}{15 (3)}  & \multicolumn{1}{l|}{}       &        \\ \hline
\end{tabular}
\end{table}

\begin{table}[!htb]
\centering
\caption{$R_{3}^{\Delta}(i,j)$ values and corresponding number of extremal graphs.}
\label{table:R_{3}(i,j)}
\begin{tabular}{ll|llllllll|}
\cline{3-10}
                                           &   & \multicolumn{8}{c|}{$j$}                                                                                                                                                                                                      \\ \cline{3-10} 
                                           &   & \multicolumn{1}{l|}{5}     & \multicolumn{1}{l|}{6}     & \multicolumn{1}{l|}{7}       & \multicolumn{1}{l|}{8}      & \multicolumn{1}{l|}{9}        & \multicolumn{1}{l|}{10}      & \multicolumn{1}{l|}{11}      & 12       \\ \hline
\multicolumn{1}{|l|}{\multirow{4}{*}{$i$}} & 5 & \multicolumn{1}{l|}{5 (7)} & \multicolumn{1}{l|}{6 (7)} & \multicolumn{1}{l|}{7 (8)}   & \multicolumn{1}{l|}{8 (8)}  & \multicolumn{1}{l|}{9 (9)}    & \multicolumn{1}{l|}{10 (9)}  & \multicolumn{1}{l|}{11 (10)} & 12 (10)  \\ \cline{2-10} 
\multicolumn{1}{|l|}{}                     & 6 & \multicolumn{1}{l|}{6 (1)} & \multicolumn{1}{l|}{7 (5)} & \multicolumn{1}{l|}{9 (1)}   & \multicolumn{1}{l|}{10 (8)} & \multicolumn{1}{l|}{12 (2)}   & \multicolumn{1}{l|}{13 (25)} & \multicolumn{1}{l|}{15 (7)}  & 16 (144) \\ \cline{2-10} 
\multicolumn{1}{|l|}{}                     & 7 & \multicolumn{1}{l|}{6 (1)} & \multicolumn{1}{l|}{8 (2)} & \multicolumn{1}{l|}{10 (1)}  & \multicolumn{1}{l|}{12 (3)} & \multicolumn{1}{l|}{15 (2)}   & \multicolumn{1}{l|}{}        & \multicolumn{1}{l|}{}        &          \\ \cline{2-10} 
\multicolumn{1}{|l|}{}                     & 8 & \multicolumn{1}{l|}{6 (1)} & \multicolumn{1}{l|}{8 (2)} & \multicolumn{1}{l|}{10 (10)} & \multicolumn{1}{l|}{13 (2)} & \multicolumn{1}{l|}{15 (551)} & \multicolumn{1}{l|}{}        & \multicolumn{1}{l|}{}        &          \\ \hline
\end{tabular}
\end{table}

\begin{table}[!htb]
\centering
\caption{$R_{4}^{\Delta}(i,j)$ values and corresponding number of extremal graphs.}
\label{table:R_{4}(i,j)}
\begin{tabular}{ll|llllllll|}
\cline{3-10}
                                           &    & \multicolumn{8}{c|}{$j$}                                                                                                                                                                                                          \\ \cline{3-10} 
                                           &    & \multicolumn{1}{l|}{6}      & \multicolumn{1}{l|}{7}      & \multicolumn{1}{l|}{8}      & \multicolumn{1}{l|}{9}       & \multicolumn{1}{l|}{10}       & \multicolumn{1}{l|}{11}        & \multicolumn{1}{l|}{12}       & 13      \\ \hline
\multicolumn{1}{|l|}{\multirow{5}{*}{$i$}} & 6  & \multicolumn{1}{l|}{6 (14)} & \multicolumn{1}{l|}{7 (14)} & \multicolumn{1}{l|}{8 (14)} & \multicolumn{1}{l|}{9 (14)}  & \multicolumn{1}{l|}{10 (14)}  & \multicolumn{1}{l|}{11 (14)}   & \multicolumn{1}{l|}{12 (14)}  & 13 (14) \\ \cline{2-10} 
\multicolumn{1}{|l|}{}                     & 7  & \multicolumn{1}{l|}{7 (1)}  & \multicolumn{1}{l|}{8 (6)}  & \multicolumn{1}{l|}{10 (1)} & \multicolumn{1}{l|}{11 (7)}  & \multicolumn{1}{l|}{12 (36)}  & \multicolumn{1}{l|}{13 (194)}  & \multicolumn{1}{l|}{14 (959)} & 16 (41) \\ \cline{2-10} 
\multicolumn{1}{|l|}{}                     & 8  & \multicolumn{1}{l|}{7 (1)}  & \multicolumn{1}{l|}{9 (2)}  & \multicolumn{1}{l|}{11 (1)} & \multicolumn{1}{l|}{12 (44)} & \multicolumn{1}{l|}{14 (20)}  & \multicolumn{1}{l|}{15 (3115)} & \multicolumn{1}{l|}{}         &         \\ \cline{2-10} 
\multicolumn{1}{|l|}{}                     & 9  & \multicolumn{1}{l|}{7 (1)}  & \multicolumn{1}{l|}{9 (2)}  & \multicolumn{1}{l|}{11 (7)} & \multicolumn{1}{l|}{13 (19)} & \multicolumn{1}{l|}{15 (146)} & \multicolumn{1}{l|}{18 (255)}  & \multicolumn{1}{l|}{}         &         \\ \cline{2-10} 
\multicolumn{1}{|l|}{}                     & 10 & \multicolumn{1}{l|}{7 (1)}  & \multicolumn{1}{l|}{9 (2)}  & \multicolumn{1}{l|}{11 (7)} & \multicolumn{1}{l|}{13 (70)} & \multicolumn{1}{l|}{16 (123)} & \multicolumn{1}{l|}{}          & \multicolumn{1}{l|}{}         &         \\ \hline
\end{tabular}
\end{table}

Lastly, we suggest the following generalization of Theorems \ref{thm:T_k(k+2)}, \ref{thm:T_k(k+3)} and \ref{thm:T_k(k+4)} as a conjecture. 
\begin{conj}\label{conj}
With the preceeding notation, $T_k(k+i)=k+2i-1$ for $2 \leq i \leq k$, with $K_{i-1,k+i-1}$ as an extremal graph.
\end{conj}
Clearly, the complete bipartite graph $K_{i-1,k+i-1}$ does not contain a $k$-sparse set of size $k+i$ for $2\leq i \leq k$. This implies $T_k(k+i)\geq k+2i-1$. For $k\geq 2$, Conjecture \ref{conj} claims all $T_k(j)$ values, where $k+2\leq j \leq 2k$. This suggests that for large (but fixed) $k$, there are $k-1$ values of $T_k(j)$ that grow linearly. Note that this trend does not continue as $n$ grows since we have $T_k(n) = \Theta(\frac{n^2}{\log n})$ for fixed $k$ by Theorem \ref{T_k(n)}.\\
 Referring to Table \ref{table:T_k(j)}, each colored diagonal corresponds to $T_k(k+i)$ values for a fixed $i$. We note that the values of $T_k(k+i)$ grow linearly as $k$ goes to infinity for fixed $i$. Moreover, the extremal graph count and their structures are the same along a diagonal (for fixed $i$). The non-colored values on said diagonals (which fall out of the range for $2 \leq i \leq k$) do not carry the observed regularity. Theorems \ref{thm:T_k(k+2)}, \ref{thm:T_k(k+3)} and \ref{thm:T_k(k+4)} prove Conjecture \ref{conj} for $i=2$ (orange), $i=3$ (yellow) and $i=4$ (blue) respectively. Furthermore, the values on the green diagonal ($i = 5$) demonstrates the regularity and supports the conjecture. We suspect that this unexpected pattern on extremal graphs continues for larger $k$. The last value we could compute is for $i=6$ is $T_6(12)=17$ with 723 extremal graphs. 


\section{Conclusions}

In the search for defective Ramsey numbers in triangle-free graphs, we have looked into two parameters which are $R^{\Delta}_k(i,j)$ and $T_k(j)$. Some defective Ramsey numbers for specific configurations of parameters ($i,j,k$) are obtained with direct proof techniques, whereas some values are obtained by computer enumeration. Further values can be developed with the aid of novel structural results for triangle-free graphs and a streamlining of our algorithms.

Growth rates of these parameters, relative to one other, is also of interest.  We do not know if, for fixed $k$, whether $T_k (m) - T_{k +1}(m)$ is bounded, let alone if the difference goes to infinity.  Similarly, we do not know the behavior of  $\frac{T_k (m)} {T_{k +1}(m)}$, for fixed $k$.  We suspect this ratio moves towards 1, but cannot prove it.  Along these same lines, we do not know if there is a small $k$ and a large $m$ where $T_k (m)= T_{k +1}(m)$. Similarly, we do not know if there is a large $m$ and small $k$ where $T_k (m)= T_k (m +1)$.

As a future work, one could investigate Conjecture \ref{conj}. This would most probably require techniques other than the one used in proving Theorems \ref{thm:T_k(k+2)}, \ref{thm:T_k(k+3)} and \ref{thm:T_k(k+4)}. Note that the number of cases for possible maximum degree values to be considered in these proofs will increase with $i$, making it inconvenient to obtain a proof for all $i$ and $k$ such that $2\leq i \leq k$ using this approach. 

In general, we think that the interaction between efficient computer enumeration methods and classical proof techniques is a promising research direction for computing defective Ramsey numbers (and/or related parameters) in various graph classes.

{\small
{\em Authors' addresses}:

{\em Tinaz Ekim}, Department of Industrial Engineering, Bo\u{g}azi\c{c}i University, 34342, Bebek, Istanbul, Turkey.
 e-mail: \texttt{tinaz.ekim@\allowbreak bogazici.edu.tr}.
 
{\em Burak Nur Erdem}, Department of Industrial Engineering, Bo\u{g}azi\c{c}i University, 34342, Bebek, Istanbul, Turkey.
 e-mail: \texttt{burak.erdem@\allowbreak bogazici.edu.tr}.

{\em John Gimbel}, Mathematics and Statistics, University of Alaska, Fairbanks, AK, 99775-6660, USA.
 e-mail: \texttt{jggimbel@\allowbreak alaska.edu}.
}

\end{document}